\documentclass[reqno,12pt] {amsart}
\usepackage{cite}

\newtheorem{theorem}{Theorem}[section]
\newtheorem{lemma}[theorem]{Lemma}
\newtheorem{prop}[theorem]{Proposition}

\theoremstyle{definition}

\theoremstyle{remark}
\newtheorem{remark}[theorem]{Remark}

\newcommand{\mysection}[1]{\section{#1}
\setcounter{equation}{0}}

\newcommand{\bR}{\mathbb R}

\newcommand\cA{\mathcal{A}}

\DeclareMathOperator{\esssup}{ess\,sup}

\renewcommand{\epsilon}{\varepsilon}

\begin{document}
\title[Partial regularity of solutions to the NSE] {Partial regularity of
solutions to the four-dimensional Navier-Stokes equations at the
first blow-up time.}

\author[H. Dong]{Hongjie Dong}
\address[H. Dong]{Department of Mathematics.
University of Chicago, 5734 S. University Avenue, Chicago, Il,
60637, USA; (current address: School of Mathematics, Institute for
Advanced Study, Einstein Drive, Princeton, NJ 08540, USA)
 } \email{hjdong@math.ias.edu}
 \thanks{The first author was partially supported by the
National Science Foundation under agreement No. DMS-0111298.}

\author[D. Du]{Dapeng Du}
\address[D. Du]{School of Mathematical Sciences,
Fudan University, Shanghai 200433, P.R. China}
\email{dpdu@fudan.edu.cn}

\date{\today}



\begin{abstract}
The solutions of incompressible Navier-Stokes equations in four
spatial dimensions are considered. We prove that the two-dimensional
Hausdorff measure of the set of singular points at the first blow-up
time is equal to zero.
\end{abstract}

\maketitle

\mysection{Introduction}

In this paper we consider both the Cauchy problem and the
initial-boundary value problem for incompressible Navier-Stokes
equations in {\em four} spatial dimensions with unit viscosity and
zero external force:
\begin{equation}
                            \label{NSeq1}
u_t+u\nabla u-\Delta u+\nabla p=0,\quad \text{div}\, u=0
\end{equation}
in a smooth domain $Q_T=\Omega\times (0,T)\subset \bR^d\times \bR$.
Boundary condition $u|_{\Omega\times [0,T]}=0$ is imposed if
$\Omega\neq \bR^d$. Here $d=4$ and the initial data $a$ is in the
closure of $\{u\in C_0^{\infty}(\Omega)\,;\,\text{div}\,u=0\}$ in
$L_d(\Omega)$ if $\Omega$ is bounded, or is in the closure of
$\{u\in C^{\infty}(\Omega)\,;\,\text{div}\,u=0\}$ in
$L_d(\Omega)\cap L_2(\Omega)$ if $\Omega=\bR^d$. The local
well-posedness of such problems is well-known (see, for example,
\cite{kato} and \cite{giga2}). The solution $u$ is locally smooth in
both spacial and time variables. We are interested in the partial
regularity of $u$ at the first blow-up time $T$.

Many authors have studied the partial regularity of solutions (in
particular, {\em weak} solutions) of the Navier-Stokes equations,
especially when $d$ is equal to three. V. Scheffer studied partial
regularity in a series of papers \cite{Sch1, Sch2, Sch4}. In three
space dimensions, he established various partial regularity results
for weak solutions satisfying the so-called local energy inequality.
For $d=3$, the notion of {\em suitable weak solutions} was first
introduced in a celebrated paper \cite{CKN} by L. Caffarelli, R.
Kohn and L. Nirenberg. They called a pair consisting of velocity $u$
and pressure $p$ a suitable weak solution if $u$ has finite energy
norm, $p$ belongs to the Lebesgue space $L_{5/4}$, $u$ and $p$ are
weak solutions to the Navier-Stokes equations and satisfy a local
energy inequality. After proving some criteria for local boundedness
of solutions, they established partial regularity of solutions and
estimated the Hausdorff dimension of the singular set. They proved
that, for any suitable weak solution $u,p$, there is an open subset
where the velocity field $u$ is H\"older continuous and they showed
that the 1-D Hausdorff measure of the complement of this subset is
equal to zero. In \cite{flin}, with zero external force, F. Lin gave
a more direct and sketched proof of Caffarelli, Kohn and Nirenberg's
result. A detailed treatment was later given by O. Ladyzhenskaya and
G. A. Seregin in \cite{OL2}. Very recently, some extended results
are obtained in Seregin \cite{Seregin} and Gustafson, Kang and Tsai
\cite{tsai}.

For $d=4$, V. Scheffer proved in \cite{Sch3} that there exists a
weak solution $u$ in  $\bR^4\times \bR^+$ such that $u$ is
continuous outside a locally closed set of $\bR^4\times \bR^+$ whose
3-D Hausdorff measure is finite. Although Scheffer's paper is not
recent, it appears to us that this is the only published results on
the partial regularity of 4-D Navier-Stokes equations.

\begin{remark}
                                        \label{remark1.1}
The weak solution considered in \cite{Sch3} doesn't verify the local
energy estimate. The existence of a weak solution satisfying the
local energy estimate is still an open problem.
\end{remark}

Now let's state our result. Instead of dealing with weak solutions,
we work on classical solutions of 4-D Navier-Stokes equations, which
are regular before they blow up. Our first result is the that the
singular set at the first blowup time is a compact set with zero 2-D
Hausdorff measure. We show this after  two partial regularity
criterions are obtained. Our proof is conceptually similar to Lin's
in \cite{flin}, but the problem is technically harder. The main
difficulty is the lack of certain compactness. We overcome it by a
novelty use of the backward heat kernel (see the proof of Lemma
\ref{lemma05.18.1}) and by the use of two appropriate scaled norms
of the pressure. It is possible because the nonlinear term is
controlled by using the Sobolev embedding theorem, although we don't
have a compact embedding here.

\begin{remark}
                                        \label{remark1.2}
In the setting of classical solutions, our result is the 4-D version
of Caffarelli, Kohn and Nirenberg's theorem in \cite{CKN}.
\end{remark}

As an application of one of the partial regularity criterions
derived in proof of the first result we get our second result: in
case $\Omega=\bR^4$ if a solution blows up, it must blow up at a
finite time.

\begin{remark}
We can prove a similar result in 3-D by using the same argument.
Detailed discussions on the long-time behavior of solutions to 3-D
Navier-Stokes can be found in J. Heywood \cite{heywood} and M.
Wiegner \cite{wiegner} (and references therein). It seems that we
need some new method to deal with the five or higher dimensional
case. To the authors' best knowledge all the existing methods on
partial regularity for the Navier-Stokes equations share the
following prerequisite condition: in the energy inequality the
nonlinear term should be controlled by the energy norm under the
Sobolev imbedding theorem.  Actually, four is the highest dimension
in which we have such condition. In five or higher dimensional, such
condition fails. Therefore, we cannot hope the existing methods work
in five or higher dimensional case.
\end{remark}

The article is organized as follows. Our main theorems (Theorem
\ref{thm3}-\ref{thm2}) are given in the following section. Some
auxiliary estimates are proved in section \ref{Sec3} and \ref{Sec4},
 among which Lemma \ref{lemma05.18.1} plays a crucial role. We give
the proof of our main theorems in the last section.

To conclude this Introduction, we explain some notation used in what
follows: $\bR^d$ is a $d$-dimensional Euclidean space with a fixed
orthonormal basis. A typical point in ${\mathbb R}^d$ is denoted by
$x=(x_1,x_2,...,x_d)$. As usual the summation convention over
repeated indices is enforced. And $x\cdot y=x_iy_i$ is the inner
product for $x,y\,\in\bR^d$. For $t>0$, we denote $H_t=\bR^4\times
(0,t)$ and  space points are denoted by $z=(x,t)$. Various constants
are denoted by $N$ in general and the expression $N=N(\cdots)$ means
that the given constant $N$ depends only on the contents of the
parentheses.

\mysection{Setting and main results}
                                        \label{Sec2}

We shall use the notation in \cite{OL2}. Let $\omega$ be a domain in
some finite-dimensional space. Denote $L_p(\omega ;\bR^n)$ and
$W^k_p(\omega;\bR^n)$ to be the usual Lebesgue and Sobolev spaces of
functions from $\omega$ into $\bR^n$. Denote the norm of the spaces
$L_p(\omega;\bR^n)$ and $W^k_p(\omega;\bR^n)$ by
$\|\cdot\|_{L_p,\omega}$ and $\|\cdot\|_{W^k_p,\omega}$
respectively. As usual, for any measurable function $u=u(x,t)$ and
any $p,q\in [1,+\infty]$,  we define
$$
\|u(x,t)\|_{L_t^pL_x^q}:=\big\|\|u(x,t)\|_{L_x^q}\big\|_{L_t^p}.
$$

For summable functions $p,u=(u_i)$ and $\tau=(\tau_{ij})$, we use
the following differential operators
$$
\partial_t u=u_t=\frac{\partial u}{\partial t},\,\,\,\,
u_{,i}=\frac{\partial u}{\partial x_i},\,\,\,\, \nabla
p=(p_{,i}),\,\,\,\,\nabla u=(u_{i,j}),
$$
$$
\text{div}\,u=u_{i,i},\,\,\,\,\text{div}\,\tau=(\tau_{ij,j}),\,\,\,\,
\Delta u=\text{div}\nabla u,
$$
which are understood in the sense of distributions.

We use the notation of spheres, balls and parabolic cylinders,
$$
S(x_0,r)=\{x\in \bR^4\vert |x-x_0|=r\},\quad S(r)=S(0,r),\quad
S=S(1);
$$
$$
B(x_0,r)=\{x\in \bR^4\vert |x-x_0|<r\},\quad B(r)=B(0,r),\quad B=B(1);
$$
$$
Q(z_0,r)=B(x_0,r)\times (t_0-r^2,t_0),\quad Q(r)=Q(0,r),\quad Q=Q(1).
$$
Also we denote means values of summable functions as follows
$$
[u]_{x_0,r}(t)=\frac{1}{|B(r)|}\int_{B(x_0,r)}u(x,t)\,dx,
$$
$$
(u)_{z_0,r}(t)=\frac{1}{|Q(r)|}\int_{Q(z_0,r)}u\,dz.
$$

In case $\Omega=\bR^d$, in a well-known paper \cite{kato} Kato
proved that the problem is locally well-posed. By known local
regularity theory for Navier-Stokes equations it can be proved that
solutions obtained by Kato's (also known as mild solutions) is
smooth in $\bR^d\times (0,T_*]$ for some $T_*>0$. Meanwhile, for
bounded $\Omega$, it is also known (see \cite{giga2}) that there
exists a unique solution $u$ of \eqref{NSeq1} satisfying
\begin{enumerate}
    \item $u\in C([0, T_*]; L_d)$, $u(0)=a$ for some $T_*>0$;
    \item $u\in C((0,T_*]; D(\cA^{\alpha}))$ for any $0<\alpha<1$;
    \item $\|\cA^{\alpha}u(t)\|=o(t^{-\alpha})$ as $t\to 0$.
\end{enumerate}
Here $\cA$ is the stokes operator. Moreover, such solution is smooth
in  $\Omega\times (0,T_*]$. In both cases, let $T=\sup T_*$ be the
first blow-up time. Then $u$ is a smooth function in $Q_T$.

Let $\eta(x)$ be a smooth function on $\bR^4$ supported in the unit
ball $B(1)$, $0\leq \eta\leq 1$ and $\eta\equiv 1$ on $\bar B(2/3)$.
Let $z_0$ be a given point in $\Omega\times (0,T]$ and $r>0$ a real
number such that $Q(z_0,r)\subset Q_T$. It's known that for a.e.
$t\in (t_0-r^2,t_0)$, in the sense of distribution one has
\begin{align*}
\Delta p&=\frac{\partial^2}{\partial x_i\partial
x_j}\big(u_iu_j\big)\\
&=\frac{\partial^2}{\partial x_i\partial
x_j}\big((u_i-[u_i]_{x_0,r})(u_j-[u_j]_{x_0,r})\big) \quad
\text{in}\,\,B(x_0,r).
\end{align*}
For these $t$, we consider the decomposition $p={\tilde
p}_{x_0,r}+h_{x_0,r}$ in $B(x_0,r)$, where $\tilde p_{x_0,r}$ is the
Newtonian potential of
$$(u_i-[u_i]_{x_0,r})(u_j-[u_j]_{x_0,r})\eta(x/r).$$ Then $h_{x_0,r}$
is harmonic in $B(x_0,r/2)$. In the sequel, we omit the indices of
$\tilde p$ and $h$ whenever there is no confusion. The following
notation will be used throughout the article:
$$
A(r)=A(r,z_0)=\esssup_{t_0-r^2\leq t\leq
t_0}\frac{1}{r^2}\int_{B(x_0,r)}|u(x,t)|^2\,dx,
$$
$$
E(r)=E(r,z_0)=\frac{1}{r^2}\int_{Q(z_0,r)}|\nabla u|^2\,dz,
$$
$$
C(r)=C(r,z_0)=\frac{1}{r^3}\int_{Q(z_0,r)}|u|^3\,dz,$$
$$
D(r)=D(r,z_0)=\frac{1}{r^3}\int_{Q(z_0,r)}|p-[h]_{x_0,r}|^{3/2}\,dz.
$$
$$
F(r)=F(r,z_0)=\frac{1}{r^2}\Big[\int_{t_0-r^2}^{t_0}
\big(\int_{B(x_0,r)}|p-[h]_{x_0,r}|^{1+\alpha}\,dx\big)
^{\frac{1}{2\alpha}}\,dt\big]^{\frac{2\alpha}{1+\alpha}},
$$
where $\alpha\in (0,1)$ is a number to be specified later. Notice
that these objects are all invariant under the natural scaling.

Here come our main results of the article.

\begin{theorem}
                    \label{thm3}
Let $\Omega$ be a smooth bounded set or the whole space $\bR^4$ and
$(u,p)$ be the solution of \eqref{NSeq1}. There is a positive number
$\epsilon_0$ satisfying the following property. Assume that for a
point $z_0\in \Omega\times {T}$  the inequality
\begin{equation}
                                    \label{eq4.1.12}
\limsup_{r\downarrow 0}E(r)\leq \epsilon_0
\end{equation}
holds. Then
$z_0$ is a regular point.
\end{theorem}

\begin{theorem}
                    \label{thm4}
Let $\Omega$ be a smooth bounded set or the whole space $\bR^4$ and
$(u,p)$ be the solution of \eqref{NSeq1}. There is a positive number
$\epsilon_0$ satisfying the following property. Assume that for a
point $z_0\in Q_{T}$ and for some $\rho_0>0$  such that
$Q(z_0,\rho_0)\subset Q_T$ and
\begin{equation}
                                    \label{eq4.1.11}
C(\rho_0)+D(\rho_0)+F(\rho_0)\leq \epsilon_0.
\end{equation}
Then $z_0$ is a regular point.
\end{theorem}

\begin{remark}
                                        \label{remark10.56}
It is worth noting that the object under estimation in condition
\eqref{eq4.1.12} involves the gradient of $u$ while the objects in
condition \eqref{eq4.1.11} involve only $u$ and $p$ themselves.
However, by using condition \eqref{eq4.1.12} one can obtain a better
estimate of the Hausdorff dimension of the set of all singular
points.
\end{remark}

\begin{theorem}
                    \label{thm1}
Let $\Omega$ be a smooth bounded set or the whole space $\bR^4$ and
$(u,p)$ be the solution of \eqref{NSeq1}. Then the 2-D Hausdorff
measure of the set of singular points in $\Omega\times {T}$ is equal
to zero.
\end{theorem}

\begin{theorem}
                    \label{thm2}
Assume $\Omega$ is the whole space $\bR^4$. Let $(u,p)$ be the
solution of \eqref{NSeq1}. If the solution does not blow up in
finite time, then $u$ is bounded and smooth in $\bR^4\times
(0,+\infty)$.
\end{theorem}

In the sequel, we shall make use of the following well-known
interpolation inequality.
\begin{lemma}
                    \label{lemmamulti}
For any functions $u\in W^1_2(\bR^4)$ and real numbers $q\in [2,4]$
and $r>0$,
$$
\int_{B_r}|u|^q\,dx\leq N(q)\big[\big(\int_{B_r}|\nabla u|^2\,
dx\big)^{q-2}\big(\int_{B_r}|u|^2\,dx\big)^{2-q/2}
$$
\begin{equation}
                        \label{eq11.30}
 +r^{-2(q-2)}\big(\int_{B_r}|u|^2\,dx\big)^{q/2}\big].
\end{equation}
\end{lemma}

Let $(u,p)$ be the solution of the Navier-Stokes  equation
\eqref{NSeq1}.
\begin{lemma}
                    \label{assum1}
(i) We have
\begin{equation}
                    \label{eq23.2.49}
u\in L_{\infty}(0,T;L_2(\Omega;\bR^4))\cap
L_2(0,T;W^1_2(\Omega;\bR^4))\cap L_3(Q_T),
\end{equation}
\begin{equation}
                    \label{eq23.2.50}
p\in L_{3/2}(Q_T).
\end{equation}

(ii) For $0<t\leq T$ and for all non-negative function
$$
\psi\in C_0^{\infty}(\Omega\times (0,\infty)),
$$ the following generalized
energy inequality is satisfied
\begin{multline}
\esssup_{0<s\leq t}\int_{\Omega}|u(x,s)|^2\psi(x,s)\,dx+2\int_{Q_t}|\nabla u|^2\psi\,dxds \\
                    \label{energy}
\leq \int_{Q_t}\{|u|^2(\psi_t+\Delta \psi)+(|u|^2+2p)u\cdot
\nabla\psi\}\,dx\,ds.
\end{multline}
\end{lemma}
{\em Sketch of the proof.} To prove
$$
u\in L_{\infty}(0,T;L_2(\Omega;\bR^4))\cap
L_2(0,T;W^1_2(\Omega;\bR^4)),
$$
 it is suffices to multiply the
first equation of \eqref{NSeq1} by $u$ and integrate by parts. By
using Lemma \ref{lemmamulti} with $q=3$ and integrating in $t$, we
obtain $u\in L_3(Q_T)$. Since in $Q_T$ it holds that
$$
\Delta p=\frac{\partial^2}{\partial x_i\partial x_j}u_iu_j,
$$
\eqref{eq23.2.50} follows from the Calder\'on-Zygmund's estimate.
Next, let's prove part (ii). First, notice that for $0<t<T$,
\eqref{energy} can be obtained by multiplying the first equation of
\eqref{NSeq1} by $u\psi$, integrating by parts and integrating with
respect to $t$. For the case when $t=T$, due to part (i) and
H\"older's inequality both sides of \eqref{energy} are finite. Then
it remains to let $t\to T$ and take the limit on both sides.

We shall prove Theorem \ref{thm1} in three steps. First, we want to
control $A,C,D,F$ in a smaller ball by the their values in a larger
ball under the assumption that $E$ is sufficiently small. Similar
results can be found in \cite{OL2} or \cite{flin} in case when the
space dimension is three.

\begin{lemma}
                    \label{lemma2.28.1}
Suppose $\gamma\in (0,1)$, $\rho>0$ are constants and
$Q(z_0,\rho)\subset Q_T$. Then we have
\begin{equation}
                                \label{eq11.35}
C(\gamma\rho)\leq N\big[\gamma^{-3}A^{1/2}(\rho)E(\rho)
+\gamma^{-9/2}A^{3/4}(\rho)E^{3/4}(\rho) +\gamma C(\rho)\big],
\end{equation}
where $N$ is a constant independent of $\gamma,\rho$ and $z_0$.
\end{lemma}

\begin{lemma}
                    \label{lemma2.28.2}
Suppose $\alpha\in (0,1/2]$, $\gamma\in (0,1/3]$, $\rho>0$ are
constants and $Q(z_0,\rho)\subset Q_T$. Then we have
\begin{equation}
                                \label{eq11.30.a}
F(\gamma\rho)\leq
N(\alpha)\big[\gamma^{-2}A^{\frac{1-\alpha}{1+\alpha}}(\rho)
E^{\frac{2\alpha}{1+\alpha}}(\rho)
+\gamma^{\frac{3-\alpha}{1+\alpha}} F(\rho)\big].
\end{equation}
where $N(\alpha)$ is a constant independent of $\gamma,\rho$ and
$z_0$. In particular, for $\alpha=1/2$ we have,
\begin{equation}
                                \label{eq11.29}
D(\gamma\rho)\leq N\big[\gamma^{-3}A^{1/2}(\rho)E(\rho)
+\gamma^{5/2} D(\rho)\big].
\end{equation}
Moreover, it holds that
\begin{equation}
                                \label{eq4.51}
D(\gamma\rho)\leq N(\alpha)\big[\gamma^{-3}(A(\rho)+E(\rho))^{3/2}
+\gamma^{(9-3\alpha)/(2+2\alpha)} F^{3/2}(\rho)\big].
\end{equation}
\end{lemma}

\begin{lemma}
                    \label{lemma2.28.3}
Suppose $\theta\in (0,1/2]$, $\rho>0$ are constants and
$Q(z_0,\rho)\subset Q_T$. Then we have
$$
A(\theta\rho)+E(\theta\rho)\leq N\theta^{-2}\big[
C^{2/3}(\rho)+C(\rho) +C^{1/3}(\rho)D^{2/3}(\rho)\big].
$$
In particular, when $\theta=1/2$ we have
\begin{equation}
                    \label{eq11.36}
A(\rho/2)+E(\rho/2)\leq N[C^{2/3}(\rho)+C(\rho)
+C^{1/3}(\rho)D^{2/3}(\rho)].
\end{equation}
\end{lemma}

As a conclusion, we obtain
\begin{prop}
                    \label{lemma05.03.02}
For any $\epsilon_0>0$, there exists $\epsilon_1>0$ small  such that
for any $z_0\in Q_T\cup (\bR^4\times \{T\})$ satisfying
\begin{equation}
                                    \label{eq3.02}
\limsup_{r\to 0}E(r)\leq \epsilon_1,
\end{equation}
we can find $\rho_0$ sufficiently small such
that
\begin{equation}
                                    \label{eq3.02.2}
A(\rho_0)+E(\rho_0)+C(\rho_0)+D(\rho_0)+F(\rho_0)\leq \epsilon_0.
\end{equation}
\end{prop}

In the second step, our goal is to estimate the values of $A,E,C$ and
$F$ in a smaller ball by the values of themselves in a
larger ball.

\begin{lemma}
                    \label{lemma05.18.1}
Suppose $\rho>0$, $\theta\in (0,1/3]$ are constants and
$Q(z_1,\rho)\subset Q_T$. Then we have
\begin{equation}
                    \label{eq18.2.31}
A(\theta\rho)+E(\theta\rho)\leq N\theta^2
A(\rho)+N\theta^{-3}[A(\rho)+E(\rho)+F(\rho)]^{3/2},
\end{equation}
where $N$ is a constant independent of $\rho,\theta$ and $z_1$.
\end{lemma}

\begin{lemma}
                    \label{lemma05.3.1}
Suppose $\rho>0$ is constants and $Q(z_1,\rho)\subset Q_T$. Then we
can find $\theta_1>0$ small such that
$$
A(\theta_1\rho)+E(\theta_1\rho)+F(\theta_1\rho)\leq
\frac{1}{2}\big[A(\rho)+E(\rho)+F(\rho)\big]
$$
\begin{equation}
                    \label{eq5.3.1}
+N(\theta_1)\big[A(\rho)+E(\rho)+F(\rho)\big]^{3/2},
\end{equation}
where $N$ is a constant independent of $\rho$ and $z_1$.
\end{lemma}

\begin{prop}
                    \label{lemma05.03.18.6}
For any $\epsilon_2>0$, there exists $\epsilon_0>0$ small such that:
if for some $z_0\in Q_T\cup (\bR^4\times \{T\})$ and $\rho_0>0$
satisfying $Q(z_0,\rho_0)\subset Q_T$ and
\begin{equation}
                                    \label{eq3.02.2.b.b}
C(\rho_0)+D(\rho_0)+F(\rho_0)\leq \epsilon_0,
\end{equation}
then for any $\rho\in (0,\rho_0/4)$ and $z_1\in Q(z_0,\rho/4)$ we
have
\begin{equation}
                                    \label{eq18.6.48}
A(\rho,z_1)+C(\rho,z_1)+E(\rho,z_1)+F(\rho,z_1)\leq \epsilon_2.
\end{equation}
\end{prop}

Finally, we apply Schoen's trick to prove the main theorems.

\mysection{Proof of Proposition \ref{lemma05.03.02}}
                                            \label{Sec3}

We will prove these lemma briefly. For more detail, we refer the
reader to \cite{OL2}.

{\bf Proof of Lemma \ref{lemma2.28.1}.} Denote $r=\gamma \rho$. We
have, by using Poincar\'e's inequality and Cauchy's inequality,
\begin{align*}
\int_{B(x_0,r)}|u|^2\,dx&=\int_{B(x_0,r)}
\big(|u|^2-[|u|^2]_{x_0,\rho}\big)\,dx
+\int_{B(x_0,r)}[|u|^2]_{x_0,\rho}\,dx\\
&\leq N\rho\int_{B(x_0,\rho)}|\nabla u||u|\,dx
+\big(\frac{r}{\rho}\big)^4\int_{B(x_0,\rho)}|u|^2\,dx\\
& \leq N\rho\Big(\int_{B(x_0,\rho)}|\nabla u|^2\,dx\Big)^{1/2}
\Big(\int_{B(x_0,\rho)}|u|^2\,dx\Big)^{1/2}\\
&\quad+\big(\frac{r}{\rho}\big)^4\int_{B(x_0,\rho)}|u|^2\,dx\\
& \leq N\rho^2A^{1/2}(\rho) \Big(\int_{B(x_0,\rho)}|\nabla
u|^2\,dx\Big)^{1/2} \\
&\quad+\big(\frac{r}{\rho}\big)^4
\Big(\int_{B(x_0,\rho)}|u|^3\,dx\Big)^{2/3}\rho^{4/3}.
\end{align*}
Owing to Lemma \ref{lemmamulti} with $q=3$ and using the inequality
above, one gets
\begin{multline*}
\int_{B(x_0,r)}|u|^3\,dx \leq N\big[\big(\int_{B(x_0,r)}|\nabla
u|^2\,dx\big) \rho A^{1/2}(\rho)\\
+\rho^3r^{-2} A^{3/4}(\rho)\big(\int_{B(x_0,\rho)}|\nabla
u|^2\,dx\big)^{3/4}
+\big(\frac{r}{\rho}\big)^4\int_{B(x_0,\rho)}|u|^3\,dx\big]
\end{multline*}

By integrating with respect to $t$ on $(t_0-r^2,t_0)$ and applying H\"older's inequality, we get
\begin{multline*}
\int_{Q(z_0,r)}|u|^3\,dz\leq N\Big[\big(\int_{Q(z_0,\rho)}|\nabla u|^2\,dz\big)\rho A^{1/2}(\rho)\\
+\rho^3 r^{-3/2} A^{3/4}(\rho) \Big(\int_{Q(z_0,\rho)}|\nabla
u|^2\,dz\Big)^{3/4}
+\big(\frac{r}{\rho}\big)^4\int_{Q(z_0,r)}|u|^3\,dz\Big].
\end{multline*}
The conclusion of Lemma \ref{lemma2.28.1} follows immediately.

{\bf Proof of Lemma \ref{lemma2.28.2}.} Denote $r=\gamma \rho$.
Recall the decomposition of $p$ introduced in Section \ref{Sec2}. By
using Calder\'on-Zygmund's estimate, Lemma \ref{lemmamulti} with
$q=2(1+\alpha)$ and Poincar\'e inequality, one has
\begin{align}
&\int_{B(x_0,r)}|{\tilde p}_{x_0,r}(x,t)|^{1+\alpha}\,dx\nonumber\\
&\leq N\int_{B(x_0,r)} |u-[u]_{x_0,r}|^{2(1+\alpha)}\,dx\nonumber\\
&\leq N\Big(\int_{B(x_0,r)}|\nabla u|^2\,dx\Big)^{2\alpha}
\Big(\int_{B(x_0,r)}|u-[u]_{x_0,r}|^2\,dx\Big)^{1-\alpha}\nonumber\\
&\quad+Nr^{-4\alpha}\Big(\int_{B(x_0,r)}
|u-[u]_{x_0,r}|^2\,dx\Big)^{1+\alpha}\nonumber \\
                    \label{eq28.1}
&\leq N\Big(\int_{B(x_0,r)}|\nabla u|^2\,dx\Big)^{2\alpha}
\Big(\int_{B(x_0,r)}|u|^2\,dx\Big)^{1-\alpha}.
\end{align}
Here we also use the inequality
$$
\int_{B(x_0,r)} |u-[u]_{x_0,r}|^2\,dx \leq \int_{B(x_0,r)}|u|^2\,dx.
$$
Similarly,
\begin{equation}
                            \label{eq065.23}
\int_{B(x_0,\rho)}|{\tilde p}_{x_0,\rho}|^{1+\alpha}\,dx \leq
N\Big(\int_{B(x_0,\rho)}|\nabla u|^2\,dx\Big)^{2\alpha}
\Big(\int_{B(x_0,\rho)}|u|^2\,dx\Big)^{1-\alpha}.
\end{equation}
Since $h_{x_0,\rho}$ is harmonic in $B(x_0,\rho/2)$, any Sobolev
norm of $h_{x_0,\rho}$ in a smaller ball can be estimated by any of
its $L_p$ norm in $B(x_0,\rho/2)$. Thus, by using poincar\'e
inequality one can obtain
\begin{align}
&\int_{B(x_0,r)}|h_{x_0,\rho}-[h_{x_0,\rho}]_{x_0,r}|^{1+\alpha}\,dx\nonumber\\
&\leq Nr^{1+\alpha}\int_{B(x_0,r)}|\nabla
h_{x_0,\rho}|^{1+\alpha}\,dx \nonumber\\
&\leq
Nr^{5+\alpha}\sup_{B(x_0,r)}|\nabla h_{x_0,\rho}|^{1+\alpha}\nonumber\\
& \leq N\big(\frac{r}{\rho}\big)^{5+\alpha}
\int_{B(x_0,\rho/2)}|h_{x_0,\rho}(x,t)-[h_{x_0,\rho}]_{x_0,\rho}|^{1
+\alpha}\,dx\nonumber\\
                    \label{eq28.2}
&\leq N\big(\frac{r}{\rho}\big)^{5+\alpha}\Big[
\int_{B(x_0,\rho)}|p(x,t)-[h_{x_0,\rho}]_{x_0,\rho}|^{1+\alpha}
+|{\tilde p}_{x_0,\rho}(x,t)|^{1+\alpha}\,dx\Big].
\end{align}
Combining \eqref{eq065.23} and \eqref{eq28.2} together yields,
\begin{align}
&\int_{B(x_0,r)}|p(x,t)-[h_{x_0,\rho}]_{x_0,r}|^{1+\alpha}\,dx\nonumber\\
&\leq N\Big(\int_{B(x_0,\rho)}|\nabla u(x,t)|^2\,dx\Big)^{2\alpha}
\Big(\int_{B(x_0,\rho)}|u(x,t)|^2\,dx\Big)^{1-\alpha}\nonumber\\
                                        \label{eq065.42}
&\quad +N\big(\frac{r}{\rho}\big)^{5+\alpha}
\int_{B(x_0,\rho)}|p(x,t)-[h_{x_0,\rho}]_{x_0,\rho}|^{1+\alpha}\,dx.
\end{align}
Since $\tilde p_{x_0,r}+h_{x_0,r}=p=\tilde
p_{x_0,\rho}+h_{x_0,\rho}$ in $B(x_0,r)$, by H\"older's inequality
\begin{align}
&\int_{B(x_0,r)}|[h_{x_0,\rho}]_{x_0,r}-[h_{x_0,r}]_{x_0,r}|^{1+\alpha}\,dx\nonumber\\
&=Nr^4|[h_{x_0,\rho}]_{x_0,r}-[h_{x_0,r}]_{x_0,r}|^{1+\alpha}\nonumber\\
&= Nr^4|[\tilde p_{x_0,\rho}]_{x_0,r}-[\tilde
p_{x_0,r}]_{x_0,r}|^{1+\alpha}\nonumber\\
                                        \label{eq065.47}
& \leq N\int_{B(x_0,r)}|{\tilde p}_{x_0,\rho}|^{1+\alpha}+|{\tilde
p}_{x_0,r}|^{1+\alpha}\,dx.
\end{align}
From \eqref{eq28.1}, \eqref{eq065.23}, \eqref{eq065.42} and
\eqref{eq065.47} we get
\begin{align}
&\int_{B(x_0,r)}|p(x,t)-[h_{x_0,r}]_{x_0,r}|^{1+\alpha}\,dx
\nonumber\\
& \leq N\Big(\int_{B(x_0,\rho)}|\nabla u(x,t)|^2\,dx\Big)^{2\alpha}
\Big(\int_{B(x_0,\rho)}|u(x,t)|^2\,dx\Big)^{1-\alpha}\nonumber\\
                                        \label{eq065.56}
&\quad +N\big(\frac{r}{\rho}\big)^{5+\alpha}
\int_{B(x_0,\rho)}|p(x,t)-[h_{x_0,\rho}]_{x_0,\rho}|^{1+\alpha}\,dx.
\end{align}
Raising to the power $1/(2\alpha)$ and integrating with respect to
$t$ in $(t_0-r^2,t_0)$ complete the proof of \eqref{eq11.30.a} and
also \eqref{eq11.29}.

To prove \eqref{eq4.51}, we use a slightly different estimate from
\eqref{eq28.2}. Again, since $h$ is harmonic in $B(x_0,\rho/2)$, we
have
\begin{align}
&\int_{B(x_0,r)}|h_{x_0,\rho}-[h_{x_0,\rho}]_{x_0,r}|^{3/2}\,dx\nonumber\\
&\leq
Nr^{3/2}\int_{B(x_0,r)}|\nabla h_{x_0,\rho}|^{3/2}\,dx\nonumber\\
&\leq Nr^{11/2}\sup_{B(x_0,r)}|\nabla h_{x_0,\rho}|^{3/2}\nonumber\\
&\leq N\frac{r^{11/2}}{\rho^{3/2+6/(1+\alpha)}}
\Big[\int_{B(x_0,\rho)}|h_{x_0,\rho/2}(x,t)-[h_{x_0,\rho}]_{x_0,\rho}|^{1+\alpha}\,dx
\Big]^{\frac{3}{2(1+\alpha)}}\nonumber\\
&\leq N\frac{r^{11/2}}{\rho^{3/2+6/(1+\alpha)}}\Big\{\Big[
\int_{B(x_0,\rho)}|p(x,t)-[h_{x_0,\rho}]_{x_0,\rho}|
^{1+\alpha}\,dx\Big]^{\frac{3}{2(1+\alpha)}}\nonumber\\
                            \label{eq18.5.07}
&\quad+\Big[\int_{B(x_0,\rho)}|{\tilde p}_{x_0,\rho}(x,t)|
^{1+\alpha}\,dx\Big]^{\frac{3}{2(1+\alpha)}}\Big\}.
\end{align}
Similar to \eqref{eq065.56}, we obtain
\begin{align}
&\int_{B(x_0,r)}|p(x,t)-[h_{x_0,r}]_{x_0,r}|^{3/2}\,dx\nonumber\\
& \leq N\Big(\int_{B(x_0,\rho)}|\nabla u(x,t)|^2\,dx\Big)
\Big(\int_{B(x_0,\rho)}|u(x,t)|^2\,dx\Big)^{1/2}\nonumber\\
&\quad+N\frac{r^{11/2}}{\rho^{3/2+6/(1+\alpha)}}
\Big\{[\int_{B(x_0,\rho)}|p(x,t)-[h]_{x_0,\rho}|^{1+\alpha}\,dx
\big]^{\frac{3}{2(1+\alpha)}}\nonumber\\
&\quad+\Big(\int_{B(x_0,\rho)}|\nabla
u(x,t)|^2\,dx\Big)^{\frac{3\alpha}{1+\alpha}}
\Big(\int_{B(x_0,\rho)}|u(x,t)|^2\,dx\Big)^{\frac{3(1-\alpha)}
{2(1+\alpha)}}\Big]\Big\}
\end{align}

Integrating with respect to $t$ in $(t_0-r^2,t_0)$ and applying
H\"older's inequality  complete the proof of \eqref{eq4.51}.

{\bf Proof of Lemma \ref{lemma2.28.3}.} Let $r=\theta \rho$. In the
energy inequality \eqref{energy}, we put $t=t_0$ and choose a
suitable smooth cut-off function $\phi$ such that
$$
\psi\equiv 0\,\,\text{in}\,Q_{t_0}\setminus Q(z_0,\rho), \quad 0\leq
\psi\leq 1\,\,\text{in}\,Q_T,
$$
$$
\psi\equiv 1 \,\,\text{in}\,Q(z_0,r),\quad |\nabla
\psi|<N\rho^{-1},\,\, |\partial_t \psi|+|\nabla^2
\psi|<N\rho^{-2}\,\,\text{in}\,Q_{t_0}.
$$
By using \eqref{energy} and because $u$ is divergence free, we get
\begin{align*}
&A(r)+2E(r)\\
& \leq \frac{N}{r^2}\Big[\frac{1}{\rho^2}
\int_{Q(z_0,\rho)}|u|^2\,dz
\frac{1}{\rho}\int_{Q(z_0,\rho)}(|u|^2+2|p-[h]_{x_0,\rho}|)|u|\,dz\Big].
\end{align*}
Due to H\"older's inequality, one can obtain
$$
\int_{Q(z_0,\rho)}|u|^2\,dz\leq
\big(\int_{Q(z_0,\rho)}|u|^3\,dz\big)^{2/3}
\big(\int_{Q(z_0,\rho)}\,dz\big)^{1/3}\leq \rho^4C^{2/3}(\rho),
$$
\begin{align*}
&\int_{Q(z_0,\rho)}|p-[h]_{x_0,\rho}||u|\,dz \\
&\leq \big(\int_{Q(z_0,\rho)}|p-[h]_{x_0,\rho}|^{3/2}\,dz\big)^{2/3}
\big(\int_{Q(z_0,\rho)}|u|^3\,dz\big)^{1/3}\\
& \leq N\rho^3D^{2/3}(\rho)C^{1/3}(\rho).
\end{align*}
Then the conclusion of Lemma \ref{lemma2.28.3} follows immediately.

{\bf Proof of Proposition \ref{lemma05.03.02}.} Let's prove first
\eqref{eq3.02.2} without the presence of $F$ on the left-hand side.
For a given point
$$
z_0=(x_0,t_0)\in Q_T\cup (\bR^4\times \{T\})$$ satisfying
\eqref{eq3.02}, choose $\rho_0>0$ such that $Q(z_0,\rho_0)\subset
Q_T$. Then for any $\rho\in (0,\rho_0]$ and $\gamma\in (0,1/6)$, by
using \eqref{eq11.36},
$$
A(\gamma\rho)+E(\gamma\rho)\leq
N[C^{2/3}(2\gamma\rho)+C(2\gamma\rho)+D(2\gamma\rho)].
$$
This estimate, \eqref{eq11.35} and \eqref{eq11.29} together with
Young's inequality imply
\begin{align}
&A(\gamma\rho)+E(\gamma\rho)+C(\gamma\rho)+D(\gamma\rho)
\nonumber\\
&\leq N[\gamma^{2/3}C^{2/3}(\rho)+\gamma^{5/2} D(\rho)+\gamma
C(\rho)+\gamma A(\rho)]\nonumber\\
&\quad+N\gamma^{-100}(E(\rho)+E^3(\rho))\nonumber\\
&\leq N\gamma^{2/3}[A(\rho)+E(\rho)+C(\rho)+D(\rho)]+N\gamma^{2/3}\nonumber\\
                                \label{eq17.11.15}
&\quad+N\gamma^{-100}(E(\rho)+E^3(\rho)).
\end{align}
It is easy to see that for any $\epsilon_3>0$, there're sufficiently
small real numbers $\gamma\leq 1/(2N)^{3/2}$ and $\epsilon_1$ such
that if \eqref{eq3.02} holds then for all small $\rho$ we have
$$
N\gamma^{2/3}+N\gamma^{-100}(E(\rho)+E^3(\rho))<\epsilon_3/2
$$
By using \eqref{eq17.11.15} we reach
$$
A(\rho_1)+C(\rho_1)+D(\rho_1)\leq \epsilon_3
$$
for some $\rho_1>0$ small enough. To include $F$ in the estimate, it
suffices to use \eqref{eq11.30.a}.

\mysection{Proof of Proposition \ref{lemma05.03.18.6}}
                                                    \label{Sec4}

{\bf Proof of Lemma \ref{lemma05.18.1}.} Let $r=\theta\rho$. Define
the backward heat kernel as
$$
\Gamma(t,x)=\frac{1}{4\pi^2(r^2+t_1-t)^2}
e^{-\frac{|x-x_1|^2}{2(r^2+t_1-t)}}.
$$
In the energy inequality \eqref{energy} we put $t=t_1$ and choose
$\psi=\Gamma\phi:=\Gamma\phi_1(x)\phi_2(t)$, where $\phi_1,\phi_2$
are suitable smooth cut-off functions satisfying
$$
\phi_1\equiv 0\,\,\text{in}\,\bR^4\setminus B(x_1,\rho),\quad 0\leq
\phi_1\leq 1\,\,\text{in}\,\bR^4,\quad \phi_1\equiv 1
\,\,\text{in}\,B(x_1,\rho/2)
$$
$$
\phi_2\equiv
0\,\,\text{in}\,(-\infty,t_1-\rho^2)\cup(t_1+\rho^2,+\infty),\quad
0\leq \phi_2\leq 1\,\,\text{in}\,\bR,
$$
$$
\phi_2\equiv 1 \,\,\text{in}\,(t_1-\rho^2/4,t_1+\rho^2/4),\quad
|\phi_2'|\leq N\rho^{-2}\,\,\text{in}\,\bR,
$$
\begin{equation}
                                \label{eq4.17}
|\nabla \phi_1|<N\rho^{-1},\,\, |\nabla^2
\phi_1|<N\rho^{-2}\,\,\text{in}\,\bR^4.
\end{equation}
By using the equality
$$
\Delta\Gamma+\Gamma_t=0,
$$
 we have
\begin{align}
&\int_{B(x_0,\rho)}|u(x,t)|^2\Gamma(t,x)\phi(x,t)\,dx
+2\int_{Q(z_0,\rho)}|\nabla u|^2\Gamma\phi\,dz \nonumber\\
&\leq \int_{Q(z_0,\rho)}\{|u|^2(\Gamma\phi_{t} +\Gamma\Delta
\phi+2\nabla
\phi\nabla \Gamma)\nonumber\\
                                \label{eq2.54}
&\quad +(|u|^2+2p)u\cdot (\Gamma\nabla\phi+\phi\nabla\Gamma)\}\,dz.
\end{align}

After some straightforward computations, it is easy to see the following three properties:\\
(i) For some constant $c>0$, on $\bar Q(z_1,r)$ it holds that
$$
\Gamma\phi=\Gamma\geq cr^{-4}.
$$
(ii) For any  $z\in Q(z_1,\rho)$, we have
$$
|\phi(z)\nabla\Gamma(z)|+|\nabla\phi(z)\Gamma(z)|\leq Nr^{-5}.
$$
(iii) For any  $z\in Q(z_1,\rho)\setminus Q(z_1,r)$, we have
$$
|\Gamma(z)\phi_t(z)|+|\Gamma(z)\Delta\phi(z)|+|\nabla\phi\nabla\Gamma|
\leq N\rho^{-6},\quad
$$
These properties together with \eqref{eq2.54} and \eqref{eq4.17}
yield
\begin{equation}
                                    \label{eq4.47}
A(r)+E(r)\leq N[\theta^2 A(\rho) +\theta^{-3}(C(\rho)+D(\rho))].
\end{equation}
Owing to Lemma \ref{lemmamulti} with $q=3$, one easily gets
\begin{equation}
                                    \label{eq4.48}
C(\rho/3)\leq NC(\rho)\leq N[A(\rho)+E(\rho)]^{3/2}.
\end{equation}
By using \eqref{eq4.51} with $\gamma=1/3$, we have
\begin{equation}
                                    \label{eq5.39}
D(\rho/3)\leq N[A(\rho)+E(\rho)+F(\rho)]^{3/2}.
\end{equation}
Upon combining \eqref{eq4.47} (with $\rho/3$ in place of $\rho$),
\eqref{eq4.48} and \eqref{eq5.39} together, the lemma is proved.

{\bf Proof of Lemma \ref{lemma05.3.1}.}  Due to \eqref{eq11.30.a} and
\eqref{eq18.2.31}, for any $\gamma,\theta\in (0,1/3]$, we have
\begin{align}
F(\gamma\theta\rho)&\leq
N\big[\gamma^{-2}(A(\theta\rho)+E(\theta\rho))
+\gamma^{(3-\alpha)/(1+\alpha)}F(\theta\rho)\big]\nonumber \\
&\leq N\gamma^{-2}\theta^2
A(\rho)+\gamma^{(3-\alpha)/(1+\alpha)}\theta^{-2}F(\rho)\nonumber \\
                                \label{eq6.18}
&\quad+N\gamma^{-2}\theta^{-3}
\big[A(\rho)+E(\rho)+F(\rho)\big]^{3/2}
\end{align}
\begin{equation}
                    \label{eq18.6.25}
A(\gamma\theta\rho)+E(\gamma\theta\rho)\leq (\gamma\theta)^2
A(\rho)+(\gamma\theta)^{-3}[A(\rho)+E(\rho)+F(\rho)]^{3/2}.
\end{equation}
Now we put $\alpha=1/27$ such that
$$
(3-\alpha)/(1+\alpha)=20/7>2.
$$
In Section \ref{Sec5}, we will give more explanation why we choose
$\alpha=1/27$. Now one can choose and fix $\gamma$ and $\theta$ sufficiently small
such that
$$
N[\gamma^{-2}\theta^2+\gamma^{20/7}\theta^{-2}+(\gamma\theta)^2]\leq
1/2.
$$
Upon adding \eqref{eq6.18} and \eqref{eq18.6.25}, we obtain
$$
A(\gamma\theta\rho)+E(\gamma\theta\rho)+F(\gamma\theta\rho)\leq
\frac{1}{2}A(\rho)+N[A(\rho)+E(\rho)+F(\rho)]^{3/2},
$$
where $N$ depends only on $\theta$ and $\gamma$. After putting
$\theta_1=\gamma\theta$, the lemma is proved.

{\bf Proof of Proposition \ref{lemma05.03.18.6}.} Take the constant
$\theta_1$ from Lemma \ref{lemma05.3.1}. Due to Lemma
\ref{lemma2.28.3}, we may choose $\epsilon_0, \epsilon'>0$ small
enough such that
$$
A(\rho_0/2)+E(\rho_0/2)+C(\rho_0/2)+D(\rho_0/2)+F(\rho_0/2)\leq
\epsilon'.
$$
\begin{equation}
                            \label{eq9.58}
2\epsilon'+8N(\theta_1)\epsilon'^{3/2} \leq
\min(4\epsilon',\theta_1^2\epsilon_2),
\end{equation}
where the constant $N(\theta_1)$ is the same one as in
\eqref{eq5.3.1}. Since $z_1\in Q(z_0,\rho/4)$, we have
$$
Q(z_1,\rho_0/4)\subset Q(z_0,\rho_0/2)\subset Q_T,
$$
$$
A(\rho_0/4,z_1)+E(\rho_0/4,z_1)+F(\rho_0/4,z_1)\leq 4\epsilon'.
$$
By using \eqref{eq9.58} and \eqref{eq5.3.1}, one obtains inductively
for $k=1,2,\cdots,$
$$
A(\theta_1^k\rho_0/4,z_1)+E(\theta_1^k\rho_0/4,z_1)
+F(\theta_1^k\rho_0/4,z_1)\leq
\min\{\theta_1^2\epsilon_2,4\epsilon'\},
$$
Thus, for any $\rho\in (0,\rho_0/4]$, it holds that
$$
A(\rho,z_1)+E(\rho,z_1)+F(\rho,z_1)\leq \epsilon_2.
$$
To include the term $C(\rho,z_1)$ in the estimate, it suffices to
use \eqref{eq4.48}. The proposition is proved.

\mysection{Proof of Theorem \ref{thm3}-\ref{thm2}}
                                \label{Sec5}
{\bf Proof of Theorem \ref{thm3} and \ref{thm4}.} Let $z_0\in
Q_T\cup (\bR^4\times \{T\})$ be a given point. Proposition
\ref{lemma05.03.02} and \ref{lemma05.03.18.6} imply that for any
$\epsilon_2>0$ there exist small numbers
$\epsilon_1,\epsilon_0,\rho_0>0$ such that either
\begin{equation}
                                        \label{eq10.30}
\limsup_{r\to 0}E(r,z_0)\leq \epsilon_1
\end{equation}
or
\begin{equation}
                                        \label{eq4.2.12}
C(\rho_0)+D(\rho_0)+F(\rho_0)\leq \epsilon_0
\end{equation}
holds true, we can find $\rho_1>0$ so that $Q(z_0,\rho_1)\subset
Q_T$ and for any $z_1\in Q(z_0,\rho_1/2),\rho\in (0,\rho_1/2)$ we
have
\begin{equation}
                                    \label{eq18.10.21}
C(\rho,z_1)+F(\rho,z_1)\leq \epsilon_2.
\end{equation}

Let $\delta\in (0,\rho_1^2/4)$ be a number and denote
$$
M_{\delta}=\max_{\bar Q(z_0,\rho_1/2)\cap \bar
Q_{T-\delta}}d(z)|u(z)|,
$$
where
$$
d(z)=\min[\text{dist}(x,\partial \Omega),(t+\rho_1^2/4-T)^{1/2}].
$$

\begin{lemma}
                            \label{lemma5.1}
If \eqref{eq18.10.21} holds true for a sufficiently small
$\epsilon_2$, then
\begin{equation}
                            \label{eq10.37}
\sup_{Q(z_0,\rho_1/4)}|u(z)|<+\infty.
\end{equation}
\end{lemma}
\begin{proof}
If for all $\delta\in (0,\rho_1^2/4)$ we have $M_\delta\leq 2$, then
there's nothing to prove. Otherwise, suppose for some $\delta$ and
$z_1\in \bar Q(z_0,\rho_1/2)\cap \bar Q_{T-\delta}$,
$$
M:=M_\delta=|u(z_1)|d(z_1)>2.
$$ Let $r_1=d(z_1)/M<d(z_1)/2$. We make the scaling as follows:
$$
\bar u(y,s)=r_1u(r_1y+x_1,r_1^2s+t_1),$$
$$
\bar p(y,s)=r_1 p(r_1 y+x_1,r_1^2 s+t_1).
$$
It's known that the pair $(\bar u,\bar p)$ satisfies the
Navier-stokes equations \eqref{NSeq1} in $Q(0,1)$. Obviously,
\begin{equation}
                            \label{eq11.30.b}
\sup_{Q(0,1)}|\bar u|\leq 2,\quad |\bar u(0,0)|=1.
\end{equation}

Due to the scaling-invariant property of our objects $A,E,C,D$ and
$F$, in what follows we look them as objects associated to $(\bar
u,\bar p)$ at the origin. For any $\rho\in (0,1]$, we have
\begin{equation}
                                    \label{eq19.10.55}
C(\rho)+F(\rho)\leq \epsilon_2.
\end{equation}
Recall what we did before in the proof of Lemma \ref{lemma2.28.2}.
Since $\bar u$ is bounded in $Q(0,1)$, we have
\begin{equation}
                            \label{eq19.11.13}
\int_{Q(0,1/3)}| \tilde {\bar p}_{0,1}|^{14}\,dz\leq
\int_{Q(0,1/2)}| \tilde {\bar p}_{0,1}|^{14}\,dz\leq N,
\end{equation}
\begin{align*}
&\int_{B(0,1/3)}|\bar h_{0,1}(z)-[\bar h_{0,1}]_{0,1/3}|^{14}\,dx\\
&\leq N\sup_{B(0,1/3)}|\nabla \bar h_{0,1}(x,t)|^{14}\\
& \leq N\Big(\int_{B(0,1/2)}|\bar h_{0,1}-[\bar
h_{0,1}]_{0,1/2}|^{28/27}\,dx\Big)^{27/2},
\end{align*}
and
\begin{align}
&\int_{Q(0,1/3)}|\bar h_{0,1}(z)-[\bar
h_{0,1}]_{0,1/3}|^{14}\,dz\nonumber \\
& \leq N\int_{-1/9}^0\Big(\int_{B(0,1/2)} |\bar h_{0,1}-[\bar
h_{0,1}]_{0,1/2}|^{28/27}\,dx\Big)^{27/2}\,dt\nonumber \\
                                    \label{eq11.22}
& \leq N(1+F^{14}(1)).
\end{align}
Estimates \eqref{eq19.11.13} and \eqref{eq11.22} yield
\begin{equation}
                            \label{eq11.26}
\int_{Q(0,1/3)}|\bar p(z)-[\bar h_{0,1}]_{0,1/3}|^{14}\,dz\leq N.
\end{equation}
Because $(\bar u,\bar p)$ satisfies the equation
$$
\bar u_t-\Delta \bar u=\text{div}(\bar u\otimes \bar u)-\nabla (\bar
p-[\bar h_{0,1}]_{0,1/3})
$$
in $Q(0,1)$. Owing to \eqref{eq11.30.b}, \eqref{eq11.26} and the
classical Sobolev space theory of parabolic equation, we have
\begin{equation}
                            \label{eq19.11.58}
\bar u\in W_{14}^{1,1/2}(Q(0,1/4)),\quad \|\bar
u\|_{W_{14}^{1,1/2}(Q(0,1/4))}\leq N.
\end{equation}
Since $1/2-6/14=1/14>0$, owing to the Sobolev embedding theorem (see
\cite{O.A.}), we obtain
$$
\bar u\in C^{1/14}(Q(0,1/5)),\quad \|\bar
u\|_{C^{1/14}(Q(0,1/5))}\leq N,
$$
where $N$ is a universal constant independent of $\epsilon_1$ and
$\epsilon_2$. Therefore, we can find $\delta_1<1/5$ independent of
$\epsilon_1, \epsilon_2$ such that
\begin{equation}
                                \label{eq19.12.23}
|\bar u(x,t)|\geq 1/2 \quad\text{in}\,\,Q(0,\delta_1).
\end{equation}
Now we choose $\epsilon_2$ small enough which makes
\eqref{eq19.12.23} and \eqref{eq19.10.55} a contradiction. The lemma
is proved.
\end{proof}
Theorem \ref{thm3} and \ref{thm4} follow immediately from Lemma
\ref{lemma5.1}.

{\bf Proof of Theorem \ref{thm1}.} Take the number $\epsilon_1$ in
Lemma \ref{lemma5.1}. Denote
$$
\Omega^*:=\{z\in \Omega\times \{T\}\,|\,\limsup_{r\downarrow
0}E(r,z)\leq \epsilon_1\}.
$$
It is well known that 2-D Hausdorff measure of $\Omega\setminus
\Omega^*$ is zero. By using Lemma \ref{lemma5.1}, for any $z\in
\Omega^*$ we can find $\rho>0$ such that $u$ is bounded in
$Q(z,\rho)$. Then there's no blow-up at $z$ and $z$ is a regular
point. The theorem is proved.

{\bf Proof of Theorem \ref{thm2}.} For any $\alpha_0\in [0,1]$,  due
to Lemma \ref{assum1} (i), the interpolation inequality
\eqref{eq11.30} with $r=+\infty, q=2(1+\alpha_0)$ and H\"older's
inequality, one can easily get
\begin{equation}
                            \label{eq19.5.31}
\|u\|_{L^t_{(1+\alpha_0)/\alpha_0}L^x_{2(1+\alpha_0)}
(\bR^4\times\bR^+)}<+\infty.
\end{equation}
Since $(u,p)$ satisfies
$$
\Delta p=\frac{\partial^2}{\partial x_i\partial x_j}(u_iu_j) \quad
\text{in}\,\,\bR^4\times \bR^+,
$$
Due to Calder\'on-Zygmund's estimate, we have
\begin{equation}
                            \label{eq19.5.32}
\|p\|_{L^t_{(1+\alpha_0)/(2\alpha_0)}L^x_{1+\alpha_0}
(\bR^4\times\bR^+)}<+\infty.
\end{equation}
Because of \eqref{eq19.5.31} with $\alpha_0=1/2$ and
\eqref{eq19.5.32} with $\alpha_0=1/2,\alpha$, and again by
Calder\'on-Zygmund's estimate, for any $\epsilon_4\in (0,1)$ we can
find $R\geq 1$ sufficiently large such that for any $z_0\in
\bR^4\times (R,+\infty)$ it holds that
\begin{equation}
                                                \label{eq11.40}
C(1,z_0)+\|p\|_{L_{3/2}^tL_{3/2}^x(Q(z_0,1))}+\|p\|
_{L^t_{(1+\alpha_0)/(2\alpha_0)}L^x_{1+\alpha_0} (Q(z_0,1))}\leq
\epsilon_4,
\end{equation}
\begin{equation}
                                                \label{eq11.40a}
\|\tilde p_{z_0,1}\|_{L_{3/2}^tL_{3/2}^x(Q(z_0,1))}+\|\tilde
p_{z_0,1}\| _{L^t_{(1+\alpha_0)/(2\alpha_0)}L^x_{1+\alpha_0}
(Q(z_0,1))}\leq \epsilon_4.
\end{equation}
Thus,
\begin{equation}
                                                \label{eq11.40b}
\|h_{z_0,1}\|_{L_{3/2}^tL_{3/2}^x(Q(z_0,1))}+\|h_{z_0,1}\|
_{L^t_{(1+\alpha_0)/(2\alpha_0)}L^x_{1+\alpha_0} (Q(z_0,1))}\leq
2\epsilon_4.
\end{equation}
After combining \eqref{eq11.40} and \eqref{eq11.40b} together, it is
clear by using H\"older's inequality that
$$
C(1,z_0)+D(1,z_0)+F(1,z_0)\leq N\epsilon_4,
$$
where $N$ is independent of $\epsilon_4$.

Then owing to Proposition \ref{lemma05.03.18.6} and Lemma
\ref{lemma5.1}, for sufficiently small $\epsilon_4$ we can find a
uniform upper bound $M_0>0$ such that for any $z_0\in \bR^4\times
(R,+\infty)$
$$
\sup_{z\in Q(z_0,1/4)}|u(z)|\leq M_0.
$$
Therefore, $u$ will not blow up as $t$ goes to infinity, and Theorem
\ref{thm2} is proved.

\section*{Acknowledgment}
The authors would like to express their sincere gratitude to Prof.
V. Sverak for pointing out this problem and giving many useful
comments for improvement. The authors would also thank to Prof. N.V.
Krylov for helpful discussions and the referee for his careful
review of the article..

\end{document}